\documentclass{amsproc}
\usepackage{graphicx,amssymb}

\copyrightinfo{2004}{Ribbon Graphs United}

\theoremstyle{definition}

\theoremstyle{remark}

\numberwithin{equation}{section}
\def\BBV           {{}^{B\!}(B^\vee)}
\def\BBVpic        {{}^{B\!}(\!B^{\!\vee}\!)}
\def\calb          {{\mathcal B}}
\def\calc          {{\mathcal C}}
\def\calca         {\mbox{${\mathcal C}_{\!A}$}}
\def\calcaa        {\mbox{${}_A\!{\,}{\mathcal C}_{\!A}$}}
\def\calh          {{\mathcal H}}
\def\calm          {{\mathcal M}}
\def\calp          {{\mathcal P}}
\def\cals          {{\mathcal S}}
\def\calv          {{\mathfrak V}}
\def\complex       {{\mathbb C}}
\def\End           {{\rm End}} 
\def\id            {\mbox{\sl id}}
\def\Hom           {{\rm Hom}} 
\def\reals         {{\mathbb R}}
\def\Rbimod        {{\mathcal B}imod(R)}
\def\Rmod          {{\mathcal M}od(R)}
\def\Rep           {{\mathcal R}ep}
\def\zet           {{\mathbb Z}}

\newcommand\includemyscaledpicture[6] {{\begin{picture}(#1,#2)(#3,#4)
                   \scalebox{.#5}{\includegraphics{#6.eps}} \end{picture}}}
\newcommand\BBB[8]{\begin{picture}(#1,#2)(#3,#4)
                     \put(10,0){\begin{picture}(0,0)
                     \includemyscaledpicture0000 {38} {2simplexBBB}\end{picture}}
                     \put(12,11){\scriptsize$ #7 $}
                     \put(30,40){\scriptsize$ #6 $}
                     \put(42,11){\scriptsize$ #5 $}
                  \end{picture}}
\newcommand\BBW[8]{\begin{picture}(#1,#2)(#3,#4)
                     \put(10,0){\begin{picture}(0,0)
                     \includemyscaledpicture0000 {38} {2simplexBBW}\end{picture}}
                     \put(12,11){\scriptsize$ #7 $}
                     \put(30,40){\scriptsize$ #5 $}
                     \put(42,11){\scriptsize$ #6 $}
                  \end{picture}}
\newcommand\WWB[8]{\begin{picture}(#1,#2)(#3,#4)
                     \put(10,0){\begin{picture}(0,0)
                     \includemyscaledpicture0000 {38} {2simplexWWB}\end{picture}}
                     \put(12,11){\scriptsize$ #7 $}
                     \put(30,40){\scriptsize$ #6 $}
                     \put(42,11){\scriptsize$ #5 $}
                  \end{picture}}
\newcommand\WWW[8]{\begin{picture}(#1,#2)(#3,#4)
                     \put(10,0){\begin{picture}(0,0)
                     \includemyscaledpicture0000 {38} {2simplexWWW}\end{picture}}
                     \put(12,11){\scriptsize$ #7 $}
                     \put(30,40){\scriptsize$ #6 $}
                     \put(42,11){\scriptsize$ #5 $}
                  \end{picture}}

\begin{document}

\title{Picard groups in rational conformal field theory}

\author{J\"urg Fr\"ohlich}
\address{Institut f\"ur theoretische Physik, ETH Z\"urich, CH--8093 Z\"urich}
\email{juerg@itp.phys.ethz.ch}

\author{J\"urgen Fuchs}
\address{Institutionen f\"or fysik, Karlstads Universitet,
Universitetsg.~5, S--65188 Karlstad}
\email{jfuchs@fuchs.tekn.kau.se}

\author{Ingo Runkel}
\address{Max-Planck-Institut f\"ur Gravitationsphysik, Am M\"uhlenberg 1,
D--14476 Golm}
\email{ingo@aei.mpg.de}

\author{Christoph Schweigert}
\address{Fachbereich Mathematik, Universit\"at Hamburg,
Bundestra\ss{}e 55, D--20146 Hamburg}
\email{schweigert@math.uni-hamburg.de}
\thanks{%
Invited talk by C.S.\ at the conference on {\em Non-commutative geometry and 
representation theory in mathematical physics\/} (Karlstad, Sweden, July 2004).
To appear in the proceedings.
\\[.2em]
J.F.\ is supported by VR under project no.\ 621--2003--2385,
and C.S.\ by the DFG project SCHW 1162/1-1.
The collaboration between J.F.\ and J.F.\ is supported in part by  grant
IG 2001-070 from STINT (Stiftelsen f\"or internationalisering av h\"ogre
utbildning och forskning).  
}

\date{November 2004}   

\begin{abstract}
Algebra and representation theory in modular tensor categories can be 
combined with tools from topological field theory to obtain a deeper 
understanding of rational conformal field theories in two dimensions:
It allows us to establish the existence of sets of consistent correlation 
functions, to demonstrate some of their properties in a model-independent 
manner, and to derive explicit expressions for OPE coefficients and coefficients
of partition functions in terms of invariants of links in three-manifolds.
\\
We show that a Morita class of (symmetric special) Frobenius algebras $A$ 
in a modular tensor category $\calc$ encodes all data needed to 
describe the correlators. A Morita-invariant formulation
is provided by module categories over $\calc$.
Together with a bimodule-valued fiber functor, the system (tensor category 
+ module category) can be described by a weak Hopf algebra. 
\\
The Picard group of the category $\calc$ can be used to construct examples of 
symmetric special Frobenius algebras. The Picard group of the category of
$A$-bimodules describes the internal symmetries of the theory and allows one
to identify generalized Kramers-Wannier dualities.
\end{abstract}

\maketitle

\section{Modular tensor categories and topological field theories}

The structure of a modular tensor category appears in a variety of
representation theoretic problems in mathematical physics. For the purposes 
of this contribution, by a modular tensor category we understand an
abelian, semi-simple $\complex$-linear tensor category $\calc$ that comes
with a braiding, i.e.\ a collection of functorial isomorphisms
  $$ c_{X,Y}^{} :\quad 
  X \,{\otimes}\, Y \xrightarrow{\cong} Y \,{\otimes}\, X \,, $$
for any pair $X,Y$ of objects of $\calc$,
and a twist, i.e.\ a collection of functorial isomorphisms
  $$ \theta_X :\quad X \xrightarrow{\cong} X $$
for any object $X$ of $\calc$, such that the following axioms hold: First,
  $$ \theta_{X\otimes Y} = c_{Y,X}^{} \circ (\theta_Y\,{\otimes}\,\theta_X)
  \circ c_{X,Y}^{} \,; $$
second, there is a compatible duality; third, there are only finitely many
isomorphism classes of simple objects -- a set of representatives of which
we denote by $\{U_i\}_{i\in I} $ -- and the tensor unit ${\bf 1}\,{=}\,U_0$
is simple. Finally,
the braiding is maximally non-degenerate in the sense that the matrix
  $$ s_{ij}^{} = {\rm tr}\, c_{U_i,U_j}^{} \circ c_{U_j,U_i}^{} $$
($i,j\,{\in}\,I$) is invertible.

Modular tensor categories arise in various contexts.
For example, the representation categories of the following algebraic 
structures can be modular tensor categories: weak quantum groups,
conformal nets of von Neumann algebras on the real line, and vertex algebras. 
In view of the role of the two latter structures in two-dimensional 
conformal quantum field theory (CFT) (see e.g.\ \cite{EVka,gago}),
it follows in particular that modular tensor categories constitute the 
axiomatization of the chiral data -- in essence, the monodromy of the
conformal blocks -- of rational CFTs.

Recently, quite a few results have been obtained that
characterize cases when representation categories are modular:

\def\leftmargini{1em}
\begin{itemize}
\item 
If $H$ is a connected $C^*$ weak Hopf algebra, then the category
of unitary representations of its double is a unitary modular tensor
category \cite{nitv}.
\item 
Similarly, the representation category of a connected ribbon factorizable 
weak Hopf algebra over $\complex$ (or, more generally, over any algebraically
closed field $\Bbbk$) with a Haar integral is modular \cite{nitv}.
\item 
If a finite-index net of von Neumann algebras on the real line is strongly 
additive (which for conformal nets is equivalent to Haag duality) and has the 
split property, its category of local sectors is a modular tensor category 
\cite{kalm}.
\item 
Finally, according to the results of \cite{huan19}, if a self-dual vertex 
algebra that obeys Zhu's $C_2$ cofiniteness condition and certain conditions
on its homogeneous subspaces has a semi-simple representation
category, then this category is actually a modular tensor category.  
\end{itemize}

The definition of a modular tensor category was
motivated \cite{tura6} by the fact 
that it allows for the construction of a three-dimensional topological
field theory (TFT). Such a TFT furnishes a modular functor, i.e.\ it
assigns finite-dimensional vector spaces to two-manifolds -- more precisely, 
to extended surfaces -- and linear maps to cobordisms.

An extended surface (for a given tensor category $\calc$) is a closed 
oriented two-dimensional manifold ${\rm X}$ with finitely many embedded (germs 
of) arcs labelled by objects of $\calc$, together with the choice of a Lagrangian 
subspace in $H_1({\rm X},\reals)$. There is a natural notion of morphisms of 
extended surfaces. Given two extended surfaces ${\rm X}$ and ${\rm Y}$, a 
cobordism ${\rm M}$ from ${\rm X}$ to ${\rm Y}$ is an oriented three-manifold 
with an embedded ribbon graph such that 
$\partial {\rm M} \,{=}\, {\rm X} \,{\sqcup}\, (-{\rm Y})$.

The complex vector spaces $\calh({\rm X})$ -- called the state spaces of the 
TFT, or the spaces of conformal blocks -- assigned to extended surfaces ${\rm X}$ 
obey $\calh(\emptyset)\,{=}\,\complex$ and the multiplicativity property
$\calh({\rm X}{\sqcup}{\rm Y})\,{=}\,\calh({\rm X})\,{\otimes_\complex}\,
\calh({\rm Y})$. A modular functor associates to each morphism 
$f{:}~{\rm X}\,{\to}\,{\rm Y}$ of extended surfaces a linear map 
$f_\sharp{:}~\calh({\rm X})\,{\to}\,\calh({\rm Y}) $, while to a cobordism 
$({\rm M},\partial_-{\rm M},\partial_+ {\rm M})$ it assigns a linear map
  $$ Z({\rm M},\partial_-{\rm M},\partial_+ {\rm M}): \quad
  \calh(\partial_-{\rm M}) \to \calh(\partial_+ {\rm M}) \,. $$
In particular, a cobordism of the form $Z({\rm M},\emptyset,\partial{\rm M})$ 
gives rise to map $\complex \,{\to}\, \calh(\partial {\rm M})$. Put differently,
a ribbon graph in a three-manifold ${\rm M}$ allows one to specify 
a vector in the space of conformal blocks on the boundary $\partial {\rm M}$.

The axioms for the linear maps $Z$ -- naturality, multiplicativity,
normalization of the cylinder, and functoriality -- have two important
consequences:

\def\leftmargini{1em}
\begin{itemize}\addtolength\itemsep{3pt}
\item Each space $\calh({\rm X})$ of conformal blocks carries a projective
      representation of the mapping class group $\mbox{\sl Map}({\rm X})$.
\item By gluing two arcs which are labelled by simple objects $U_j$ and
      $U_j{}_{}^{\!\vee}$ via a tube with embedded $U_j$-ribbon,
      one obtains isomorphisms $\bigoplus_{j\in I}\calh({\rm X}_j)\,
      {\xrightarrow{\scriptscriptstyle\cong}}\,\calh({\rm X}')$,
      called {\em factorization rules\/}.
\end{itemize}

To formulate a TFT one employs a cobordism category of topological manifolds. 
The use of the term {\em conformal\/} block therefore needs to be justified.
Given a conformal vertex algebra $\calv$ and an $m$-tuple 
$(\calh_{\lambda_1},\calh_{\lambda_2},...\,,\calh_{\lambda_m})$ of 
$\calv$-modules, conformal blocks are constructed \cite{FRbe} as vector bundles 
$\calb(\lambda_1,\lambda_2,...\,,\lambda_m)$ with a projectively flat connection 
over the moduli space $\calm_{g,m}$ of complex curves of genus $g$ with 
$m$ marked points. The monodromy data of conformal blocks on $\complex{\rm P}^1$
can then be used to equip the representation category $\Rep(\calv)$ with
additional structure like a tensor product (fusion) and a braiding.
In certain cases, e.g.\ the ones described in \cite{huan19}, this endows 
$\Rep(\calv)$ with the structure of a modular tensor category. 
From this modular tensor category, one obtains
representations of the mapping class groups. It is an important
open conjecture that these representations are isomorphic to the
ones provided by the monodromies of the vector bundles $\calb$.
We will assume that this conjecture holds true; this allows us to pass
tacitly between topological categories of topological two-manifolds
and categories of conformal or complex manifolds.

The vector bundles of conformal blocks are -- up to the choice of local 
coordinates, a subtlety we ignore for the purposes of this contribution -- 
constructed as subbundles of the trivial bundle 
  $$\calm_{g,m}\times 
  (\calh_{\lambda_1}\otimes_\complex^{} \calh_{\lambda_2} \otimes_\complex^{}
  \cdots \otimes_\complex^{} \calh_{\lambda_m})^* . $$
Applying a flat section of the subbundle to an $m$-tuple of vectors
$v_i\,{\in}\,\calh_{\lambda_i}$ therefore yields a multivalued function on 
$\calm_{g,m}$. A central question, to be addressed in the next section,
is how these multivalued functions
are related to physical correlation functions of the conformal field theory.


\section{Geometry for correlators}\label{sec:geom}

To describe the correlators of a local conformal field theory on a surface 
${\rm X}$ -- which may have a non-empty boundary and, for the moment, is 
supposed to be oriented -- we associate to ${\rm X}$ its (oriented) double 
$\hat{\rm X}$.  The double $\hat {\rm X}$ comes with an orientation reversing 
involution $\sigma$ such that ${\rm X}=\hat {\rm X} / \{1,\sigma\}$. If we work 
in a geometric category, i.e.\ if ${\rm X}$ is supposed to be a conformal 
manifold, then the double has a complex structure and $\sigma$ is anti-conformal.

Correlators on ${\rm X}$ are specific vectors $C({\rm X})$ in the spaces
$\calh(\hat {\rm X})$ (this is
known as the `principle of holomorphic factorization'; it has been derived 
\cite{witt39} for important classes of rational CFTs like (gauged) WZW models 
by introducing a background gauge field on the world sheet). 
These particular elements of $\calh(\hat {\rm X})$ satify two types 
of conditions:

\def\leftmargini{1em}
\begin{itemize}\addtolength\itemsep{2pt}
\item 
The vector $C({\rm X})\,{\in}\,\calh(\hat {\rm X})$ is required to be invariant 
under the 
mapping class group $\mbox{\sl Map}({\rm X})$, which can be identified with the 
subgroup of $\mbox{\sl Map}(\hat {\rm X})$ that commutes with $\sigma$. (Recall 
that $\mbox{\sl Map}(\hat{\rm X})$ acts projectively on $\calh(\hat{\rm X})$.)
\item 
The correlators must obey factorization constraints. We do not write down these
explicitly here; in short, the image of a correlator under an isomorphism that 
is given by the factorization rules for the conformal blocks on the double is 
again a correlator. (See \cite{fjfrs} for precise statements and proofs.)
\end{itemize}

In view of the properties of TFT discussed above, it is natural to look for an 
oriented three-manifold ${\rm M}_{\rm X}$ with boundary $\partial{\rm M}_{\rm X}
\,{=}\,\hat {\rm X}$ such that $C({\rm X})\,{:=}\,
Z({\rm M}_{\rm X},\emptyset,\hat {\rm X})(1)$ is the correlator. 
The following ansatz turns out to be successful: take the product of the double 
$\hat {\rm X}$ with the interval $[-1,1]$ and mod out by the orientation 
preserving $\zet_2$ action for which the non-trivial element acts like $\sigma$ 
on $\hat {\rm X}$ and like $t\,{\mapsto}\, {-}t$ on the interval:
  $$ {\rm M}_{\rm X} = \left( \hat {\rm X} \times [-1,1] \right)
  / \,(\sigma,t\,{\mapsto}\, {-}t) \, . $$
This so-called connecting manifold \cite{fffs2} has boundary 
$\partial{\rm M}_{\rm X}\,{=}\,\hat{\rm X}$; there is a natural embedding 
$\iota{:}~{\rm X}\,{\to}\, {\rm M}_{\rm X}$ acting as $x\,{\mapsto}\, [x,0]$, 
which shows that ${\rm X}$ is a retract of ${\rm M}_{\rm X}$.  
For example, when ${\rm X}$ is a disk, ${\rm M}_{\rm X}$ is a full ball.

To specify a correlator $C({\rm X})\,{\in}\, \calh(\hat {\rm X})$, the 
three-manifold ${\rm M}_{\rm X}$ must be endowed with a ribbon graph. 
Points in the interior of ${\rm X}$ have two pre-images in $\hat {\rm X}$,
points on $\partial {\rm X}$ only one. For each field insertion point $p$ 
we place a ribbon along the distinguished interval in ${\rm M}_{\rm X}$ that 
joins $\iota(p)$ to the pre-image(s) of $p$; we also put ribbons along the 
boundary of $\iota({\rm X})$. All ribbons are labelled by objects of $\calc$, 
in a manner to be described in the next section.
Finally, we triangulate $\iota({\rm X})$, in such a manner that only trivalent
vertices occur. (Faces can have an arbitrary number of edges; thus, properly 
speaking, we are dealing with the dual of a triangulation.)
Each marked point $\iota(p)$ for $p$ in the interior of ${\rm X}$ must lie on 
an edge of the triangulation $T_{\rm X}$. We place ribbons along the edges of 
$T_{\rm X}$. This gives in particular rise to trivalent and quadrivalent 
intersections, at which we put coupons on which the ribbons end 
-- trivalent for vertices of the triangulation and for the points 
$\iota(p)$ with $p \,{\in}\, \partial {\rm X}$, and quadrivalent for the points 
$\iota(p)$ with $p$ in the interior of ${\rm X}$. 

For each of the coupons an element in the appropriate morphism space 
must be chosen. For the part of the triangulation that lies in the
interior of  $\iota({\rm X})$ this is done as follows. Choose an object
$A$ of $\calc$.  On each edge we place two ribbons labeled by $A$ which
start at one of the vertices and have orientation pointing away from that
vertex, and join them with a suitable morphism  $\Phi$ in
$\Hom(A,A^\vee)$, to be specified below. At each vertex we need an
element in $\Hom({\bf 1},A\,{\otimes}\,A\,{\otimes}\,A)$. A comparison of
this situation with the analysis of so-called lattice TFTs
\cite{fuhk,bape}, where  $\calc$ is the category of finite-dimensional
$\complex$-vector spaces, leads us to expect that $A$ is a generalization
of a Frobenius algebra to  more general tensor categories $\calc$. As we
will see, this allows us to  give a model-independent approach to
correlators of rational CFTs that is  based on a combination of TFT in
three dimensions and of non-commutative algebra in tensor categories.


\section{Frobenius algebras in modular tensor categories}\label{sec:frob}

It is in fact not hard to generalize many notions of algebra and
representation theory from vector spaces (or modules over commutative rings)
to more general tensor categories. A (unital, associative) algebra 
$A \equiv (A,m,\eta)$ in a tensor
category $\calc$, for example, is an object $A$ of $\calc$ together
with a multiplication morphism $m\,{\in}\,\Hom(A\,{\otimes}\,A,A)$ that obeys
  $$ m\circ (\id_A\otimes m) = m\circ (m\otimes \id_A)$$
and together with a unit morphism $\eta\,{\in}\,\Hom({\bf 1},A)$ such that
  $$ m\circ (\eta\otimes \id_A) = \id_A = m\circ (\id_A \otimes \eta) \, . $$
A (coassociative, counital) coalgebra $(A,\Delta,\varepsilon)$ is defined 
analogously. A {\em Frobenius algebra\/} in $\calc$ is an algebra that is also 
a coalgebra, with the additional property that the coproduct
$\Delta\,{\in}\,\Hom(A,A\,{\otimes}\,A)$ is a morphism of $A$-bimodules.

The Frobenius algebras of interest to us possess two more properties.
First, as in most applications in representation theory, they are 
{\em symmetric\/}: denote by $d_A\,{\in}\,\Hom({\bf 1},A\,{\otimes}\,A^\vee)$ and 
$\tilde d_A\,{\in}\,\Hom({\bf 1},A^\vee\,{\otimes}\,A)$ the two coevaluations 
of the category $\calc$ (which we assume to be sovereign). There are
two natural isomorphisms -- in fact isomorphisms of $A$-bimodules -- 
  $$ \Phi_1:= \left( (\varepsilon\,{\circ}\,m)\,{\otimes}\,\id_{A^\vee} \right)
  \circ ( \id_A\,{\otimes}\,d_A ) ,\qquad
  \Phi_2 := \left( \id_{A^\vee}\,{\otimes}\,(\varepsilon\,{\circ}\,m) \right)
  \circ ( \tilde d_A\,{\otimes}\,\id_A )  $$
in $\Hom(A,A^\vee)$; in a symmetric Frobenius algebra, these two isomorphisms 
coincide. It is the morphism $\Phi \equiv
\Phi_1\,{=}\,\Phi_2$ that we use along the edges of the triangulation
of $\iota({\rm X})$. Second, our Frobenius algebras are {\em special\/}, 
which means that $\Delta$ is a right-inverse of the multiplication -- this 
generalizes the notion of a separable algebra over a field -- and that 
$\varepsilon\,{\circ}\,\eta\,{=}\,\dim(A)\,\id_{\bf 1}$. 

\medskip

It can be shown \cite{fuRs4} that in a rational CFT the operator product 
(OPE) for boundary fields that preserve a given boundary condition leads 
to a symmetric special Frobenius algebra in the modular tensor category
$\calc$ that describes the chiral data of the CFT. The main ingredients are 
the associativity of the OPE and the non-degeneracy of the two-point 
correlation function of two boundary fields on a disk, which furnishes the 
non-degenerate invariant bilinear form. It should also be appreciated that a 
Frobenius algebra obtained this way from boundary fields in CFT
is not necessarily (braided-)commutative, and that
the underlying boundary condition is not required to be `elementary'.

With this in mind, the main insight of our construction can be summarized as 
follows \cite{fuRs}. For given chiral data $\calc$, a full local CFT -- which
we denote as CFT($A$) -- can be constructed from a symmetric special Frobenius 
algebra $A$ in $\calc$. This Frobenius algebra is the algebra of 
boundary fields (or, in string theory terminology, of open string states)
for a single boundary condition.

As for other boundary conditions than the one used to define $A$, the analysis 
of boundary OPEs involving also boundary condition changing operators shows 
\cite{fuRs4} that they correspond to modules over the Frobenius algebra $A$.
Here modules are defined in the obvious way: a (left) $A$-module is a pair
$M \,{\equiv} (\dot M,\rho_M)$ consisting of an object $\dot M$ of $\calc$ and 
a morphism $\rho_M\,{\in}\,\Hom(A\,{\otimes}\,\dot M,\dot M)$ such that
$\rho_M\,{\circ}\,(m\,{\otimes}\,\id_{\dot M}) \,{=}\, \rho_M\,{\circ}\,
(\id_A\,{\otimes}\,\rho_M)$ and $\rho_M\,{\circ}\,(\eta\,{\otimes}\,\id_{\dot M})
\,{=}\, \id_{\dot M}$. Many aspects of representation theory can be generalized,
see e.g.\ \cite{pare,kios,fuSc16} (in fact, some peculiar aspects can only be 
seen in a braided setting, compare \cite{ffrs}). 
For instance, there is a notion of induced module, reciprocity theorems hold, 
every simple module appears in the decomposition of an induced module, and one
can show that the module category of a special Frobenius algebra in a 
semi-simple tensor category is again semi-simple.

\medskip

These observations supply us with the first few entries in a dictionary 
relating physical concepts to algebraic notions: boundary conditions are 
$A$-modules, `elementary' boundary conditions are simple $A$-modules; 
a direct sum of simple $A$-modules amounts to introducing boundary conditions 
with Chan-Paton multiplicities.  For our construction of a ribbon graph in 
the connecting manifold ${\rm M}_{\rm X}$, we conclude that ribbons labelled 
with the object underlying a boundary condition are to be placed along the 
boundary segments of $\iota({\rm X})$.

This dictionary can be extended, and this extension at the same time
completes our labeling of the ribbons in the connecting manifold
${\rm M}_{\rm X}$. A boundary field $\Psi^{MN}_U$ that changes the
boundary condition from $M$ to $N$ has a single chiral insertion $U$ and a 
trivalent vertex at the boundary of $\iota({\rm X})$ that must be labeled by 
an element of $\Hom_A(M\,{\otimes}\,U,N)$, where $M\,{\otimes}\,U$
carries the obvious structure of a left $A$-module and $\Hom_A$ denotes
morphisms of left $A$-modules. 

Field insertions $p$ in the interior of ${\rm X}$ have
two pre-images. Bulk fields are thus labeled $\Phi_{U V}$: the
two ribbons in ${\rm M}_{\rm X}$ that originate from the pre-images of $p$ 
are inward-pointing and are labeled by $U$ and $V$, respectively. 
Further, these ribbons hit the $A$-ribbon that is placed on the 
edge of the triangulation $T_{\rm X}$ passing through $\iota(p)$, and
the corresponding quadri-valent vertex is labeled by an element of
$\Hom_{\!A|A}(U\,{\otimes}^+ A\,{\otimes^-}\,V,A)$, a morphism of $A$-bimodules.
Here the superscripts $\pm$ indicate that the object $U\,{\otimes}\,A\,{\otimes}
\,V$ is given the following structure of an $A$-bimodule: the left action is
$(\id_{U}\,{\otimes}\,m\,{\otimes}\,\id_{V})\circ ( c^{-1}_{U,A} 
    \,{\otimes}\,\id_{A\otimes V})$, while the right $A$-action is
$(\id_{U}\,{\otimes}\,m\,{\otimes}\,\id_{V})\circ (\id_{U\otimes A} 
    \,{\otimes}\,c^{-1}_{A,V} )$.

It is natural to not only consider the special $A$-bimodule
$A$ itself, but allow for arbitrary $A$-bimodules $B$ as well. They correspond 
to (tensionless) conformal defect lines which can be added to the
triangulation. Defect fields can change such defects; the
corresponding quadri-valent vertices for a change of defect from
$B_1$ to $B_2$ are labeled by $A$-bimodule morphisms in
$\Hom_{A|A}(U\,{\otimes}^+\,B_1\,{\otimes^-}\,V,B_2)$.

Using the ansatz for obtaining the correlators in terms of
three-manifolds with embedded ribbon graphs described above,
factorization and invariance under the action of the mapping class group 
$\mbox{\sl Map}({\rm X})$ can be proven \cite{fjfrs}. Also, other mathematical
objects defined previously in the discussion of rational conformal field 
theory, like the so-called classifying algebra or NIMreps, are recovered 
naturally and their properties can be established rigorously \cite{fuRs4}. 
It is worth emphasizing that our formalism provides a unified treatment of 
all modular invariant torus partition functions -- both those of simple 
current type and exceptional modular invariants.

Note that we have started from a single arbitrary boundary condition to 
construct the (symmetric special) Frobenius algebra $A$.
A different boundary condition will, in general, give us a different
Frobenius algebra $A'$. But as it turns out, for a given CFT any two such
Frobenius algebras are Morita equivalent and, moreover, Morita equivalent
Frobenius algebras give equivalent correlators; we express this as
CFT($A'$)\,$\cong$\,CFT($A$).

Finally, we mention that our construction can be extended \cite{fuRs8} to
the situation that ${\rm X}$ is not oriented, and possibly not even orientable.
In that case, the Frobenius algebra $A$ must be equipped with the additional 
structure of a {\em Jandl algebra\/}. For a Jandl algebra $A$
there is an isomorphism of algebras $\sigma{:}~A^{\rm opp}\,{\to}\, A$ which
squares to the twist, $\sigma^2\,{=}\,\theta_A$. In the special case that 
$\calc$ is the category of vector spaces, the structure of a Jandl algebra 
reduces to a (symmetric special Frobenius) algebra with an involution.


\section{Relation to other structures}

Next we wish to describe the relation of the approach to rational CFT 
based on symmetric special Frobenius algebras to other algebraic
structures whose relevance to the problem has been suggested.

As a first step, we notice that the category \calca\ of left modules
over an algebra $A$ in a tensor category $\calc$ carries the structure
of a so-called module category \cite{pare,beRn} over $\calc$: there
is a ``mixed'' tensor functor
  $$ \otimes:\quad \calca \times \calc \to \calca  $$
with an associativity constraint that satisfies generalized triangle and
pentagon axioms. Morita equivalent algebras can be characterized
\cite{ostr} by the fact that they yield equivalent module categories.

For a symmetric special Frobenius algebra in an abelian semi-simple tensor 
category, the module category $\calm\,{=}\,\calca$ is abelian and semi-simple. 
We can therefore find a complex semi-simple
algebra $R$ whose representation category $\Rmod$ is equivalent, as an
abelian category, to $\calm$. Obviously, $R$ is a direct sum of
full matrix algebras whose number of minimal ideals equals the number
of simple objects in $\calm$. Since only the number of minimal ideals, but 
not their dimension, matters, $R$ is a non-canonical object. In any case, 
this equivalence endows $\Rmod$ with the structure of a module category over 
$\calc$.

Next, recall the elementary fact that an abelian group $M$ is a module over a 
ring $S$ iff there is a morphism of rings from $S$ into ${\rm End}(M)$.
An analogous theorem is valid for categories \cite{ostr}; the ring
$S$ is replaced by the tensor category $\calc$, the module $M$ by the
module category $\calm$, and morphisms of rings by fiber functors, i.e., in the 
setting we are interested in,
by monoidal functors. If $\calm\,{\cong}\,\Rmod$ 
as an abelian category, then the bimodules $\Rbimod$ play the role of 
${\rm End}(M)$, and we have a natural bijection between fiber functors from 
$\calc$ into $\Rbimod$ and structures of a module category over $\calc$ on 
$\Rmod$.

Thus in the situation of interest to us
we obtain a bimodule-valued fiber functor
  $$ \omega_R: \quad \calc \to \Rbimod \,. $$
One would now like to apply familiar arguments from Tannaka theory to the 
algebra $H_R\,{:=}\,{\rm End}(\omega_R)$ of endomorphisms of the fiber functor 
to endow it with some structure that generalizes Hopf algebras. This can indeed
be done \cite{szla5}, provided that separability data for $R$ are chosen, 
i.e.\ a right-inverse of the multiplication and a left-inverse of the unit.  
(Except for the case that all minimal ideals of $R$ are one-dimensional,
there is no canonical separability structure.) It turns out that
$H_R$ can then be endowed with the structure of a weak Hopf algebra and
that $\calm od(H_R)$ is equivalent, as a tensor category, to $\calc$.

This construction has a converse: any weak Hopf algebra $H$ gives rise
to the tensor category $\calc\,{:=}\,{\calm od}(H)$ of left $H$-modules and
a module category $\calm$ over $\calc$: A weak Hopf algebra contains
two commuting associative unital subalgebras $H_s$ and $H_t$ which are
related by the antipode. Since the antipode is an anti-morphism of
algebras, $H_t$ can be identified with the opposite algebra of $H_s$.
By restriction to $H_s$ and $H_t$, any left $H$-module can be seen to be an
$H_t$-bimodule, and hence we have found a tensor functor from $\calc$ to 
$\calb imod(H_t)$. The general result mentioned above now implies that
the category ${\calm od}(H_t)$ is a module category over $\calc$.

Now indeed, as has been argued in \cite{bppz,pezu6}, the structure of
a weak Hopf algebra (called Ocneanu double triangle algebra) 
can be abstracted from a rational CFT and its boundary conditions. 
The discussion above shows that this algebra is not canonical; there are
infinitely many non-isomorphic weak Hopf algebras which lead to equivalent
module categories and hence to one and the same CFT.
It is also worth mentioning that in this description the braiding on $\calc$ 
must be expressed in terms of an $R$-matrix for the weak Hopf algebra, a 
structure that is not as well amenable to explicit computations as our pictorial
calculus using ribbon graphs. Still, this approach can provide non-trivial 
insight; for instance, using the fact that Davydov--Yetter cohomology
of the pair $\calm,\calc$ can be expressed in terms of Hochschild cohomology 
of any of the Hopf algebras $H_R$, it was shown in \cite{etno} that
rational conformal field theories cannot be deformed within the
class of rational conformal field theories.

\medskip

Given a module category $\calm\,{\cong}\,\calca$ over a tensor category $\calc$, 
immediately a third category enters the game the tensor category $\calc^*\,
{\cong}\,\calcaa$ of $A$-bimodules.  In our situation the latter category is 
actually equivalent to the category of module functors $\calm\,{\to}\,\calm$ 
and thus does not depend on the choice of $A$ in a Morita class. 
In contrast to $\calm$ the category $\calc^*$ is a tensor category,
with tensor unit $A$. While in our case $\calc$ is braided, $\calc^*$
is not braided, in general. Indeed, the objects of $\calc^*$
have the physical interpretation of tensionless defect lines 
whose fusion \cite{pezu5,pezu6} is described by the tensor product on $\calc^*$,
and there is no reason for the fusion of defects to be commutative.

In the present situation we then deal with four bifunctors \cite{ostr}: 
the tensor products of $\calc$ and $\calc^*$, the defining operation 
$\calm\,{\times}\, \calc\,{\to}\, \calm$ of the module category $\calm$, 
and a functor
  $$ \calc^* \times \calm \to \calm $$
acting on objects as $(B,M) \,{\mapsto}\, B\,{\otimes_{\!A}^{}}\,M$. There 
are five associativity constraints: one for each of the tensor categories
$\calc$ and $\calc^*$, and three mixed constraints for the threefold
products $\calm\,{\times}\,\calc\,{\times}\,\calc$, 
$\calc^*\,{\times}\,\calc^*\,{\times}\,\calm$
(stating that $\calm$ is a module category over both
$\calc$ and $\calc^*$) and $\calc^*\,{\times}\,\calm\,{\times}\,\calc$. 
Associativity of higher products is ensured by six
axioms of pentagon type \cite[Section 4.3]{ostr}. It should be emphasized
that the braiding on $\calc$ -- which is a crucial structure in the 
application to CFT -- is not accounted for in this setup.

\medskip

Other notions that have been discussed in this context are graphs and cells 
\cite{pezu6}. The term `cell' is motivated by the following visualization of the
structure described above in terms of oriented simplices in three dimensions. 
Vertices can be coloured ``black'' or ``white''. This gives three types of 
edges: those joining two black vertices are to be labelled by simple objects of 
$\calc$, those joining two white vertices by simple objects of $\calc^*$. 
An edge joining a white vertex to a black one is labelled by a simple object 
of $\calm$. (In principle, one can also admit edges joining a black vertex
to a white vertex. If the module category $\calm$ is realized as a category
of left $A$-modules, then one should label these edges by right $A$-modules.)

These edges can form triangles; one only needs to consider triangles for which the edges
are oriented in such a way that the boundary does not form a closed oriented 
path. 
The corresponding composition of two objects to a third object is naturally
interpreted as a (possibly mixed) tensor product.
Concretely, we have the following possibilities:
            \newpage  
\begin{tabular}{lll} {~}\\[-1.9em]
\BBB{55}{39}0{20} U V {\!W} \alpha
    \hspace{0.6em} & $ \in\ \Hom_{\calc}(U\,{\otimes}\,V,W)$
    \hspace{1.9em} & ``chiral coupling'' \\[1.02em]
\BBW{55}{39}0{20} U M N \beta
                 & $ \in\ \Hom_{\calm}(M\,{\otimes}\,U,N)$
                 & ``boundary field'' \\[1.02em]
\WWB{55}{39}0{20} B M N \gamma
                 & $ \in\ \Hom_{\calm}(B\,{\otimes}_{\!A}^{}\,M,N)$
                 & ``fusion of a defect line $B$
                   \\&&~\,to the boundary''\\[.02em]
\WWW{55}{39}0{20} B C D \delta
                 & $ \in\ \Hom_{\calc^*_{}}(B\,{\otimes}_{\!A}^{}\,C,D)$
                 & ``coupling of defect lines'' \\
\end{tabular}
\\[2.9em]
The triangles are to be labelled with elements of a basis of the morphism space
that is given in the table. The interpretation of these morphism spaces
in conformal field theory is also indicated in the table.
As is familiar from Ponzano--Regge calculus, $6j$-symbols with respect to
these bases are described in terms of scalars associated to tetrahedra
whose faces are triangles of the types shown above. 
According to the labelling of the vertices, there are then five types of 
tetrahedra. In this description, all the six pentagon axioms are interpreted
as follows: Glue two tetrahedra 
along a common face, and cut them again along an additional edge that connects
the two vertices not belonging to that face, which results in three tetrahedra.

These observations are conceptually clarified when using the language of 
two-categories (see e.g.\ section 4 of \cite{lawr3}
for an introduction). The relevant two-category has two objects $\bullet$ and 
$\circ$, corresponding to the black and white vertex. 
The morphism sets of a two-category are categories; the
endomorphism sets, in particular, are tensor categories. Thus the tensor 
category $\calc$ can be identified with the endomorphism set of $\bullet$,
and the tensor category $\calc^*$ with the endomorphism set of $\circ$.
For $x,y \,{\in}\, \{\bullet,\circ\}$, the
category $\Hom(x,y)$ is naturally a left-module category over 
the tensor category $\End(x)$ and a right-module category over $\End(y)$. 
Moreover, the category $\Hom(\circ,\bullet)$ is just $\calm$, while
$\Hom(\bullet,\circ)$ is the category of right $A$-modules. 

The whole situation can be understood as a category-theoretic analogue of a 
Morita context (compare \cite{muge8}). In fact, applying the $K_0$ functor gives
us a Morita context of complex algebras: the Grothendieck groups of the tensor 
categories $\calc$ and $\calc^*$ are even rings, $K_0(\calc)$ being
a commutative ring. $K_0(\calm)$ is a right module
over the ring $K_0(\calc)$ and a left module over the ring $K_0(\calc^*)$.
Thus it plays the role of an interpolating bimodule.

The problem of finding an ``inverse'' of the $K_0$ functor
for a given algebraic structure is known as 
{\em categorification\/}. In the case at hand, it is equivalent to the problem
of finding consistent values for all tetrahedra. In many respects, 
categorification seems to behave like a cohomology theory. Indeed, as described 
in the next section, on the Picard subcategory of $\calc$ it reduces to 
questions about the cohomology of abelian groups.


\section{Picard groups}

One might worry at this point whether there are interesting examples
of symmetric special Frobenius algebras. In every modular tensor category,
the tensor unit $\bf 1$ provides such an example; the corresponding full CFT
is also known as the ``Cardy case''. A larger class of examples 
is provided by the following general construction.

In a tensor category with duality, the isomorphism classes of invertible 
objects, i.e.\ of objects such that $V\,{\otimes}\,V^\vee{\cong}\, {\bf 1}$, 
form a group, the Picard group ${\rm Pic}(\calc)$. We denote a set of 
representatives of these isomorphism classes
     by $\{L_g \,|\, g\,{\in}\,{\rm Pic}(\calc)\}$, with $L_e\,{=}\,{\bf 1}$.
If the category $\calc$ is braided, then the Picard group is abelian.
In the physics literature, the invertible objects are known as
{\em simple currents\/} \cite{scya}.

Technically, it is convenient to consider the full subcategory 
$\calp ic(\calc)$ of direct sums of invertible objects in $\calc$ .
The Grothendieck group of this Picard subcategory is just the group
ring of the Picard group,
  $$ K_0(\calp ic(\calc)) \cong \zet {\rm Pic}(\calc) \, . $$
In this situation, categorification amounts to the following task:
given a group $G$, find a category $\calc$, such that $K_0(\calc)\,{=}\,\zet G$.
This problem has a close cousin that is of independent interest:
given an {\em abelian\/} group $G$, find a {\em braided\/} tensor
category such that $K_0(\calc)\,{=}\,\zet G$.
As it turns out, both problems possess nice answers in terms of suitable
cohomology theories. Inequivalent categorifications of a group $G$
correspond to elements of $H^3(G,\complex^\times)$ in group cohomology, while 
inequivalent braided categorifications of an abelian group $G$ are described by 
Eilenberg and Mac Lane's abelian group cohomology $H^3_{{\rm ab}}
(G,\complex^\times)$ (see \cite{fuRs9} for more explanation and references). 

It is an important result that elements of $H^3_{{\rm ab}}(G,\complex^\times)$
are in natural bijection to quadratic forms on $G$ with values in
$\complex^\times$. In fact, the braided tensor structure of the Picard 
category $\calp ic(\calc)$ is characterized by the twist of $\calc$, 
which gives a quadratic form 
$g\,{\mapsto}\,\theta_g \,{\equiv}\, \theta_{L_g}$ on ${\rm Pic}(\calc)$.
(The value $\theta_{gh} / (\theta_{g}\theta_{h})$  of the associated bilinear
form on ${\rm Pic}(\calc)$ is called the (exponentiated)
{\em monodromy charge\/} of $L_g$ with respect to $L_h$.)
As a consequence, for Picard categories the chiral 
data are particularly well accessible; this is one of the sources of the power 
of simple current methods in CFT (see e.g.\ \cite{scya,scya6,fusS6,bant6}).

We are now in a position to construct non-trivial examples of symmetric special
Frobenius algebras. We call such an algebra {\em simple\/} iff it is
simple as a bimodule over itself. (For the associated full CFT, simplicity
amounts to the property that the CFT has a unique vacuum.) A stronger condition 
is that such an algebra is simple as a left module over itself; in that case we 
call the algebra {\em haploid\/}. Now from Frobenius-Perron theory, one can 
derive the estimate $\dim_\complex {\rm Hom}(U,A) \,{\leq}\, \dim (U) $.  
for haploid Frobenius algebras.  This estimate is particularly stringent for 
those algebras for which any simple subobject $U$ is invertible and hence has 
$\dim(U)\,{=}\,1$: in these algebras the multiplicity of any simple subobject 
is either 0 or 1.  We call 
such a haploid symmetric special Frobenius algebra a {\em Schellekens algebra}.

The associativity constraint of $\calc$ gives a three-cochain 
$\psi$ on ${\rm Pic}(\calc)$. All Schel\-le\-kens algebras in a modular tensor 
category $\calc$ can be constructed by finding a subgroup $H$ of the Picard 
group ${\rm Pic}(\calc)$ and a two-cochain
$\omega{:}~H\,{\times}\, H\,{\to}\, \complex^\times$ with the property that
${\rm d} \omega \,{=}\, \psi_{|H}$. It turns out that this can 
be done if and only if for every $h\,{\in}\,H$ the twist obeys 
$(\theta_h)^{N_h}_{}\,{=}\,1$, where $N_h$ is the order of $h$.

Since the two-cochain $\omega$ depends on various gauge choices, it is somewhat 
awkward to work with $\omega$. It is therefore helpful to remember that for an
{\em abelian} group $G$ the second cohomology group $H^2(G,\complex^\times)$
-- which classifies twisted group algebras -- is canonically isomorphic to the 
group $AB(G,\complex^\times)$ of alternating bihomomorphisms on $G$ with values 
in $\complex^\times$. The isomorphism sends a representative $\omega$ of a 
cohomology class to its commutator two-cocyle $\xi$, which is defined by 
  $\xi(g,h)\,{:=}\, {\omega(g,h)}/{\omega(h,g)}$.

In the braided setting we are interested in, the notion of alternating 
bihomomorphism must be generalized; the relevant generalization is the notion 
of a {\em Kreuzer--Schellekens bihomomorphism\/} (KSB) $\Xi$, obeying
  $$ \Xi(g,h) \; \Xi(h,g) = \frac{\theta_{g}\, \theta_{h}}
  {\theta_{{gh}}} \, . $$
Here the right hand side is not equal to 1 any more, but rather is expressed in 
terms of the twist, i.e.\ of the quadratic form that characterizes the
Picard category.

The crucial observation is now that the multiplication on a Schellekens
algebra $A$ supplies us with a KSB $\Xi_A$ on the support of $A$, i.e.\ on 
the subgroup
$H(A)\,{:=}\, \{ g\,{\in}\,G \,|\, \dim_\complex\Hom(L_g,A)\,{>}\,0 \}$
of ${\rm Pic}(\calc)$, via the following relation which we display graphically:
\\ \mbox{$\ $\hspace{11.1em}}
  \begin{picture}(12,79)
  \put(1.5,0.0)  {\scalebox{.28}{\includegraphics{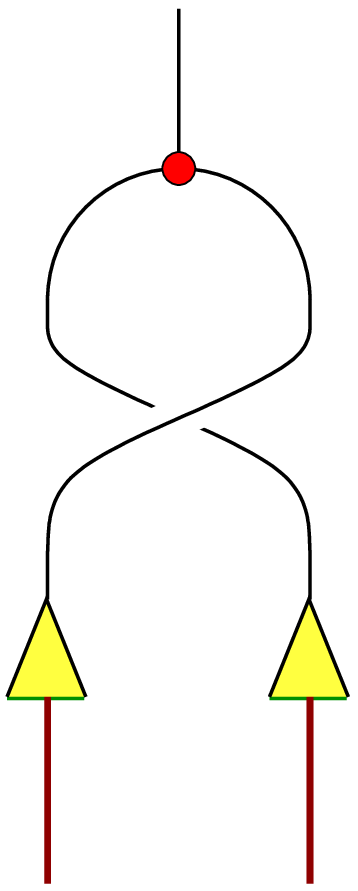}}}
  \put(111.5,0.0){\scalebox{.28}{\includegraphics{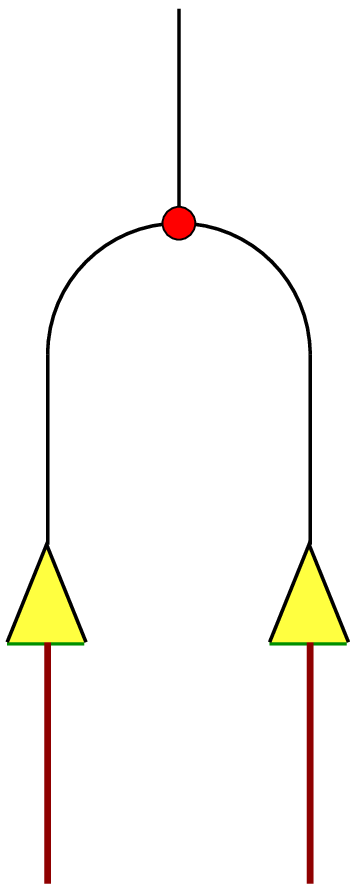}}}
      \put(0,0){\begin{picture}(0,0)
  \put(1.5,-8.5)   {\scriptsize$ L_g $}
  \put(12.9,74.5)  {\scriptsize$ A $}
  \put(17.5,60.3)  {\scriptsize$ m $}
  \put(22.7,-8.5)  {\scriptsize$ L_h $}
      \end{picture}}
  \put(48,31)      {$=:\; \Xi_A(h,g) $}
      \put(110,0){\begin{picture}(0,0)
  \put(1.5,-8.5)   {\scriptsize$ L_g $}
  \put(12.9,74.5)  {\scriptsize$ A $}
  \put(17.5,56.3)  {\scriptsize$ m $}
  \put(22.7,-8.5)  {\scriptsize$ L_h $}
      \end{picture}}
  \end{picture}

\vskip 1.4em
\noindent
Here $m$ is the multiplication morphism of $A$, and the triangles 
indicate non-zero morphisms from simple objects into $A$. The two graphs
are thus morphisms in the one-dimensional space 
$\Hom(L_g\,{\otimes}\, L_h,A)$ and hence proportional.

Conversely, one can show \cite{fuRs9} that a Schellekens algebra is uniquely
characterized by its support $H$ and a KSB on the group $H$. Hence, 
Schellekens algebras are a generalization of twisted group algebras to the 
braided setting.

It is now a central goal to express as many quantities of a local
CFT built from a Schellekens algebra as possible in terms of the KSB and other
computable quantities. This way, one obtains proofs for various simple current 
formulae that had been conjectured in the literature.

One example is the Kreuzer--Schellekens formula \cite{krSc} for the torus 
partition function, which reads
  $$ Z_{ij}(A) = \frac1{|H(A)|} \sum_{g,h\in H(A)}
  \chi_{U_i}^{}(h) \cdot \Xi_A(h,g) \cdot \delta_{\bar \jmath,gi} \,, $$
where for each simple object $U$ of $\calc$, $\chi_U^{}$ -- also called the 
monodromy charge of $U$ -- is the character of $H(A)$ given by
$\chi_U^{}(g)\,{:=}\,\theta_{L_g\otimes U} / ( \theta_{g} \theta_{U} )$.
(Note that for any simple $U$ and any $g\,{\in}\,{\rm Pic}(\calc)$,
$L_g\,{\otimes}\, U$ is again a simple object. Thus the Picard group 
${\rm Pic}(\calc)$ acts on the set of isomorphism classes of simple objects 
of $\calc$.)

Elementary boundary conditions are simple $A$-modules, which are obtained 
from the decomposition of induced $A$-modules. One finds that a simple 
$A$-module $M \,{\equiv}\, M_{O_U,\psi}$
is described by an orbit $O_U$ of the action of ${\rm Pic}(\calc)$ on the 
isomorphism classes of simple objects and an irreducible representation $\psi$ 
of the twisted group algebra $\complex_{\epsilon_U} \cals_U$. Here
$\cals_U$ is the stabilizer of the action of ${\rm Pic}(\calc)$ on the
simple object $U$, and the twist of the group algebra is characterized by the 
alternating bihomomorphism $\epsilon_U(g,h)\,{:=}\,\Phi_U(g,h)\cdot\Xi_A(h,g)$, 
where $\Phi_U$ is a certain gauge independent $6j$-symbol. This way, one
reproduces the results of \cite{fhssw} for the labelling of boundary conditions.
Similar formulae can be derived \cite{fuRs9} for defect lines, boundary states 
and other quantities in the theory.

\medskip

It is natural to consider the Picard group ${\rm Pic}(\calcaa)$ of 
invertible bimodules as well. This group has a nice physical interpretation 
\cite{ffrs3}: the corresponding defects $B$
act by internal symmetries on the correlators of the theory; explicitly:

\def\leftmargini{2.2em}
\begin{itemize}\addtolength\itemsep{2pt}
\item[(i)\phantom{ii}] 
A boundary condition described by a left module $M$ is mapped
to the one described by the left module 
  $$ B\,{\otimes_{\!A}}\, M \,{=:}\,^{B\!}M . $$
\item[(ii)\phantom{i}] 
A boundary field $\Psi^{M_1,M_2}_U$ is mapped to a boundary
field $\Psi^{^{B\!}M_1,^{B\!}M_2}_U$, where the degeneracy spaces are related
by the obvious maps $\Hom(M_1\,{\otimes}\,U, M_2) $ $
{\to}\,\Hom((B\,{\otimes_{\!A}}\,M_1)\,{\otimes}\,U, B\,{\otimes_{\!A}}\,M_2)$.
  \\[-.7em]
\item[(iii)] 
The action on bulk fields is given by 
\\ \mbox{$\ $\hspace{2.9em}}
  \begin{picture}(12,105)
  \put(0,0)      {\scalebox{.33}{\includegraphics{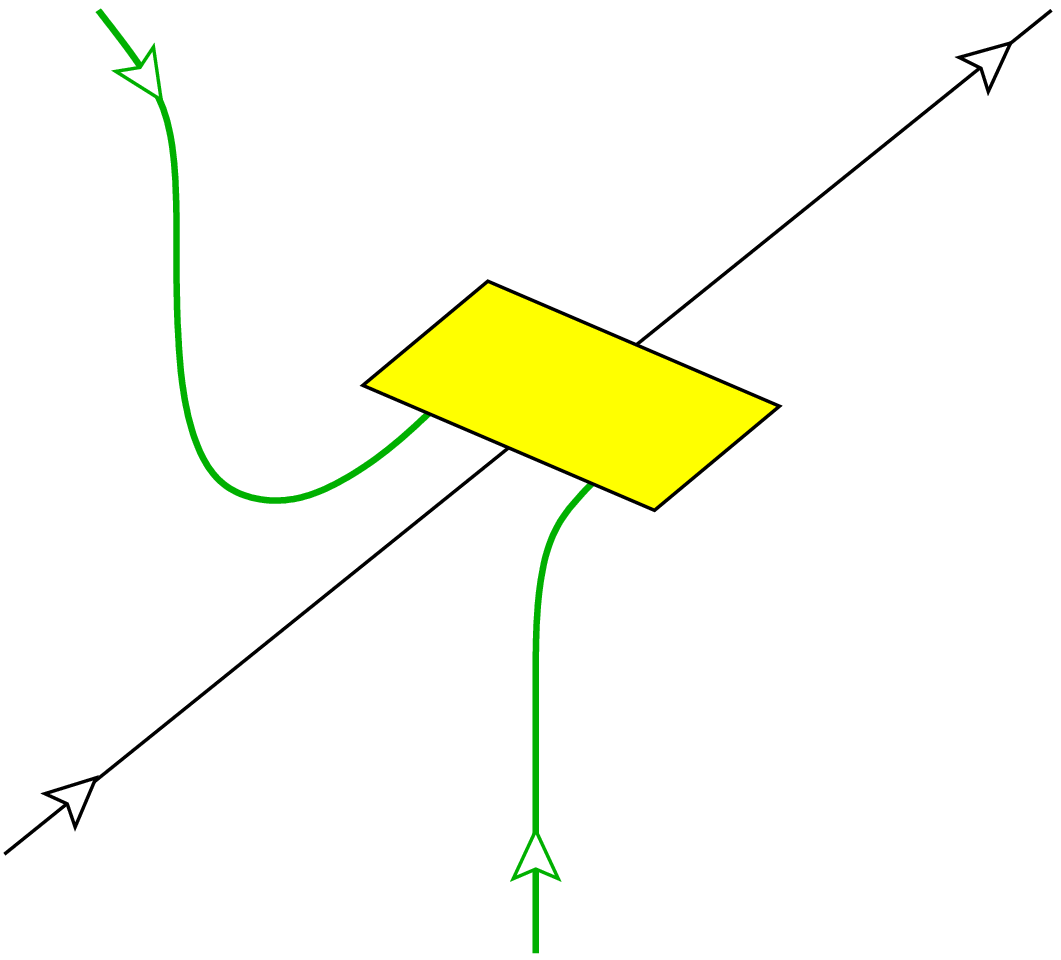}}}
  \put(0,0)      {\setlength{\unitlength}{.33pt}
                  \put( 24, 30){\scriptsize$ A $} 
                  \put(271,269){\scriptsize$ A $} 
                  \put( 49,263){\scriptsize$ U $} 
                  \put(164,  9){\scriptsize$ V $} 
                  \put(153,158){\scriptsize$ \phi $} 
                  \setlength{\unitlength}{1pt}}
  \put(153,0)    {\scalebox{.33}{\includegraphics{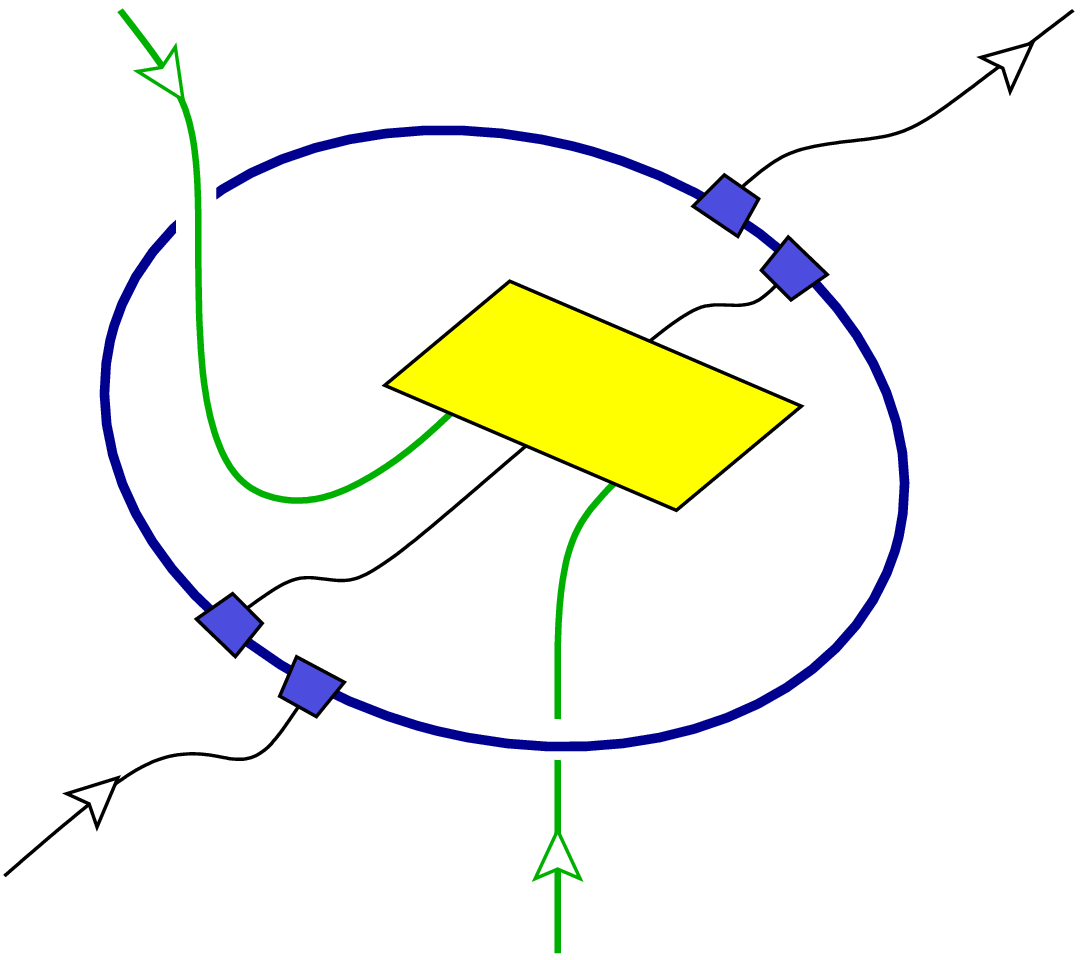}}}
  \put(153,0)    {\setlength{\unitlength}{.33pt}
                  \put( 26, 30){\scriptsize$ A $} 
                  \put(271,269){\scriptsize$ A $} 
                  \put( 50,267){\scriptsize$ U $} 
                  \put(168,  9){\scriptsize$ V $} 
                  \put(160,158){\scriptsize$ \phi $} 
                  \put(255, 83){\scriptsize$ B $} 
                  \put(244,197){\scriptsize$ \rho_l $} 
                  \put( 45, 80){\scriptsize$ \rho_l $} 
                  \put(195,235){\scriptsize$ \rho_r $} 
                  \put( 88, 57){\scriptsize$ \rho_r $} 
                  \setlength{\unitlength}{1pt}}
  \put(117,40.0) {$ \longmapsto $} 
  \end{picture}

\vskip .7em
\noindent
which defines an endomorphism of the vector space $\Hom_{A|A}(U\,{\otimes}^+A\,
{\otimes^-}\,V,A)$. Here the morphisms $\rho_{l/r}$ denote the left and right 
action of $A$ on $B$, respectively.
\end{itemize}

\noindent
With the ansatz for CFT correlators described in sections \ref{sec:geom} and 
\ref{sec:frob},
it is easy to check that this action preserves the correlation functions:
without changing the value of the correlator, one can insert in $\iota({\rm X})$
an unknot ribbon labeled by $B$. Since
$B^\vee{\otimes_{\!A}}\,B \,{\cong}\, A$, one can use contour deformation 
arguments familiar from complex analysis. This way,  $B$-loops will run 
parallel to each boundary component and encircle bulk insertions.
Everything is still connected by a network of $A$-ribbons, so that all
tensor products are to be taken over $A$. This gives precisely the 
transformation rules presented above.

It is therefore appropriate to identify elements of ${\rm Pic}(\calcaa)$
with internal symmetries of the theory. This is confirmed by the computation
of this Picard group for concrete models: for the critical Ising model one
obtains  $\zet_2$, with the non-trivial element corresponding to a global 
flip of the Ising spin, while
for the critical 3-state Potts model one obtains the symmetric group $S_3$
that permutes the three possible values of the Potts spin.


\section{Order-disorder duality from bimodules}

The discussion above shows that any bimodule in ${\rm Pic}(\calcaa)$ 
describes a symmetry of CFT$(A)$ that takes the form of equalities between 
different correlators of CFT$(A)$. Let us see how this construction gets 
modified if instead we take a defect labelled by an arbitrary 
bimodule $B\,{\in}\,{\mathcal O}bj(\calcaa)$. As it turns out, even in this 
much more general situation we get equalities between (sums of) correlators. 

To see this, we start from the correlator for a given world sheet ${\rm X}$ 
and insert in $\iota({\rm X})$ a small annular ribbon labelled by $B$. This 
changes the correlator by a factor ${\rm dim}^\calc(\dot B)/{\rm dim}^\calc(A)$,
which is the dimension of $B$ as an object of $\calcaa$.
The $B$-loop divides ${\rm X}$ in regions ``outside'' and ``inside'' the loop.
Now deform the loop until the ``outside'' area has shrunk to zero
(here we assume ${\rm X}$ to be connected). It is not difficult to check that
by this procedure of ``applying the defect $B$ to the world sheet ${\rm X}$''
one recovers the action on boundary conditions and boundary fields
described in (i) and (ii) above. 

However, in addition there are two new
effects. First, if ${\rm X}$ has a non-contractible cycle, additional defect 
lines labelled by $B{\otimes_{\!A}}B^\vee\,{\equiv}\,\BBV$  appear. For
instance, for ${\rm X}$ an annulus with boundary conditions $M$ and $N$ one gets 
\\ \mbox{$\ $\hspace{4.1em}}
  \begin{picture}(12,62)
  \put(1.5,0.0)   {\scalebox{.33}{\includegraphics{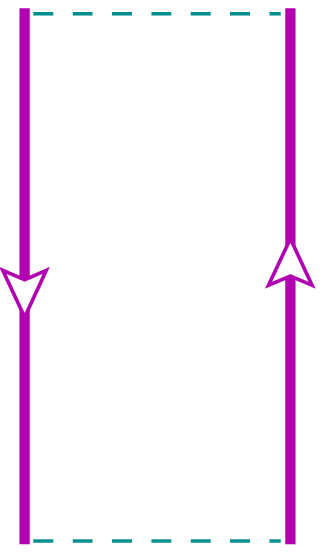}}}
  \put(77.5,0.0)  {\scalebox{.33}{\includegraphics{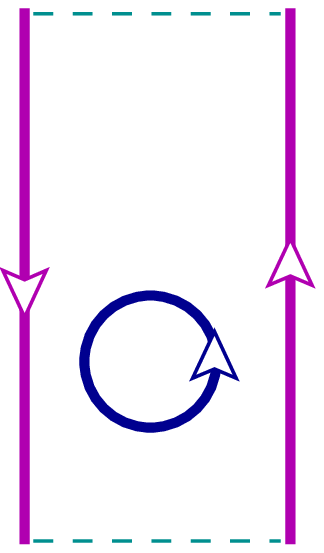}}}
  \put(153.5,-3.0){\scalebox{.33}{\includegraphics{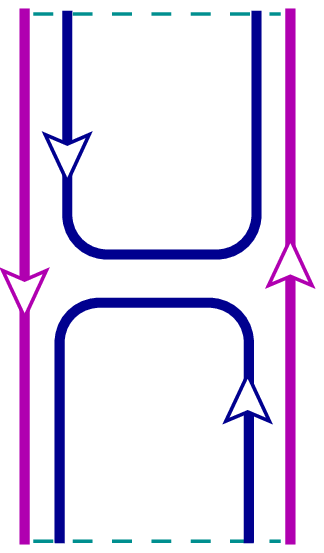}}}
  \put(229.5,0.0) {\scalebox{.33}{\includegraphics{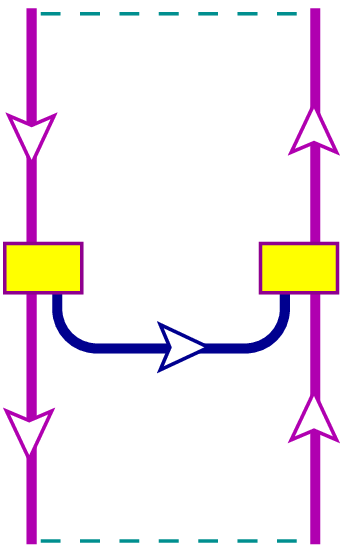}}}
      \put(1.5,0){\begin{picture}(0,0)
  \put(-7.6,33)    {\scriptsize$ M $}
  \put(29.6,17)    {\scriptsize$ N $}
      \end{picture}}
  \put(45,24)      {$ =~~~\displaystyle\frac1d $}
      \put(77.5,0){\begin{picture}(0,0)
  \put(-7.6,33)    {\scriptsize$ M $}
  \put(15,26)      {\scriptsize$ B $}
  \put(29.6,17)    {\scriptsize$ N $}
      \end{picture}}
  \put(121,24)     {$ =~~~\displaystyle\frac1d $}
      \put(153.5,0){\begin{picture}(0,0)
  \put(-7.6,33)    {\scriptsize$ M $}
  \put(29.6,17)    {\scriptsize$ N $}
  \put(13,17)      {\scriptsize$ B $}
      \end{picture}}
  \put(197,24)     {$ =~~~\displaystyle\frac1d $}
      \put(229.5,0){\begin{picture}(0,0)
  \put(-12,43)    {\scriptsize$ ^{B\!}M $}
  \put(33,43)    {\scriptsize$ ^{B\!}N $}
  \put( 7,10)      {\scriptsize$ \BBVpic $}
  \put(-5,26)    {\scriptsize$ \alpha $}
  \put(33,25)    {\scriptsize$ \beta $}
      \end{picture}}
  \end{picture}

\vskip .9em
\noindent
with $d\,{=}\,{\rm dim}^\calc(\dot B)/{\rm dim}^\calc(A)$ and suitable module 
morphisms $\alpha\,{\in}\,\Hom_{\!A}(^{B\!}M , \BBV  
   $\linebreak[0]$
{\otimes_{\!A}}{}^{B\!}M )$ and
$\beta\,{\in}\,\Hom_{\!A}(\BBV{\otimes_{\!A}}\,^{B\!}N,{}^{B\!\!}N)$.
Second, when taking the defect past a bulk field $\Phi$, in general one turns 
the bulk field into a {\em disorder field\/} $\Theta$, according to
\\ \mbox{$\ $\hspace{11.8em}}
  \begin{picture}(12,57)
  \put(1.5,0.0)   {\scalebox{.33}{\includegraphics{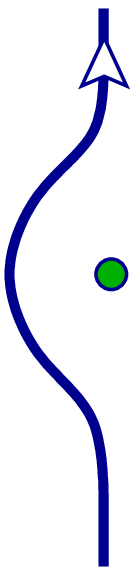}}}
  \put(77.5,0.0)  {\scalebox{.33}{\includegraphics{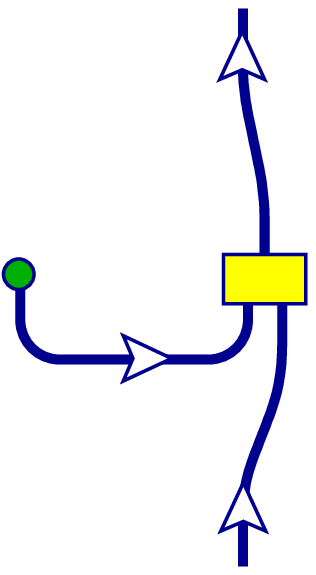}}}
      \put(0,0){\begin{picture}(0,0)
  \put(-4.6,33)   {\scriptsize$ B $}
  \put(13.8,31)   {\scriptsize$ \Phi_{UV} $}
      \end{picture}}
  \put(44,24)      {$ = $}
      \put(76,0){\begin{picture}(0,0)
  \put(-6.6,33.1) {\scriptsize$ \Theta_{UV} $}
  \put(28,10)     {\scriptsize$ B $}
  \put(28,40)     {\scriptsize$ B $}
  \put(32,27)     {\scriptsize$ \alpha $}
  \put(2,12)      {\scriptsize$ \BBVpic $}
      \end{picture}}
  \end{picture}

\vskip .9em
\noindent
with $\Theta_{UV}\,{\in}\,\Hom_{\!A|A}(U\,{\otimes}^+_{}A\,{\otimes^-}_{}\,V,B)$
a morphism of $A$-bimodules describing a disorder field, and 
$\alpha\,{\in}\,\Hom_{\!A|A}( \BBV {\otimes_{\!A}}\,B,B)$ a coupling of defect 
lines. Applying a generic defect $B$ to a world sheet thus yields an equality 
between a correlator of bulk fields and a correlator of disorder fields.

For this relationship to constitute an actual order-disorder duality symmetry, 
we must also be able to turn the disorder correlator back into
a correlator of bulk fields. In other words, we need the existence
of another defect $B'$ such that first taking $B$ past some bulk field
$\Phi$ and afterwards taking $B'$ past the resulting disorder field
$\Theta$ results in a sum of bulk field $\Phi^\beta$; pictorially,
\\ \mbox{$\ $\hspace{6.2em}}
  \begin{picture}(12,62)
  \put(1.5,0.0)   {\scalebox{.33}{\includegraphics{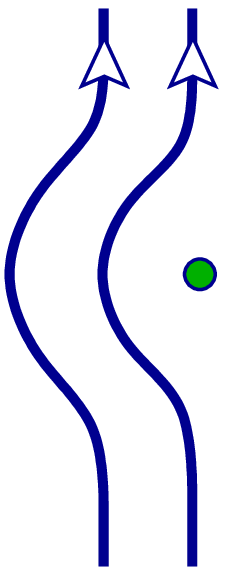}}}
  \put(84.5,0.0)  {\scalebox{.33}{\includegraphics{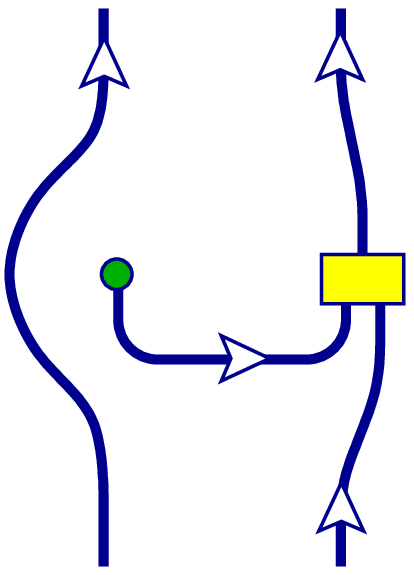}}}
  \put(183.5,0.0) {\scalebox{.33}{\includegraphics{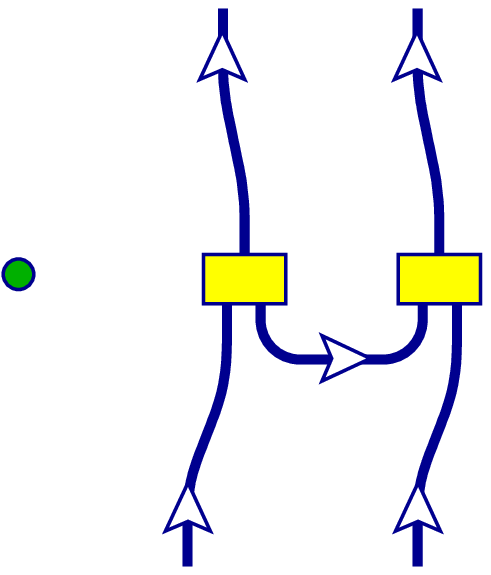}}}
      \put(0,0){\begin{picture}(0,0)
  \put(0.4,8)      {\scriptsize$ B' $}
  \put(20.5,31)    {\scriptsize$ \Phi_{UV} $}
  \put(21.1,8)     {\scriptsize$ B $}
      \end{picture}} 
  \put(54,24)      {$ = $}
      \put(83,0){\begin{picture}(0,0)
  \put(0.9,8)      {\scriptsize$ B' $}
  \put(8.5,32.1)   {\scriptsize$ \Theta_{UV} $}
  \put(35.8,8)     {\scriptsize$ B $}
  \put(42,27)  {\scriptsize$ \alpha $}
  \put(13,12)  {\scriptsize$ \BBVpic $}
      \end{picture}}
  \put(144,24)     {$\displaystyle = ~ \sum_\beta $}
      \put(182,0){\begin{picture}(0,0)
  \put(-6.1,34)  {\scriptsize$ \Phi_{UV}^\beta $}
  \put(12.1,12)    {\scriptsize$ B' $}
  \put(44.4,12)    {\scriptsize$ B $}
  \put(14,27)  {\scriptsize$ \beta $}
  \put(48,27)  {\scriptsize$ \alpha $}
      \end{picture}}
  \end{picture}
\vskip .9em
\noindent
One can prove that such a bimodule $B'$ exists iff 
  $$ B^\vee{\otimes_{\!A}}B \in {\mathcal O}bj(\calp ic(\calcaa)) \,,
  $$
in which case one can take $B'\,{=}\,B^\vee$. This generalizes the
condition $B^\vee {\otimes_{\!A}} B \,{\cong}\, A$ that must hold for 
the symmetry generating bimodules discussed in the previous section.
We call a defect line labelled by a bimodule $B$ obeying
$B^\vee{\otimes_{\!A}}B\,{\in}\,{\mathcal O}bj(\calp ic(\calcaa))$
a {\em duality defect\/}. Applying a duality defect to a world sheet results 
in a Kramers--Wan\-nier like duality, and indeed one
reproduces known dualities in this way \cite{ffrs3}.

The concept of a duality defect can still be generalized further. 
Consider two simple symmetric special Frobenius algebras $A_1$ and $A_2$.
Objects $B$ of the category ${}_{A_1\!}{\mathcal C}_{\!A_2}$ of
$A_1$-$A_2$-bimodules label tensionless conformal interfaces between 
${\rm CFT}(A_1)$ and ${\rm CFT}(A_2)$ \cite{fuRs4}. Furthermore, one checks 
that if $B^\vee{\otimes_{\!A_1}}B \,{\in}\,
{\mathcal O}bj(\calp ic({}_{A_2}\!{\,}{\mathcal C}_{\!A_2}))$ and
$B{\otimes_{\!A_2}}B^\vee \,{\in}\,
{\mathcal O}bj(\calp ic({}_{A_1}\!{\,}{\mathcal C}_{\!A_1}))$,
then $B$ gives rise to an order-disorder duality symmetry as above, 
equating this time an order correlator (i.e.\ involving no defect fields) of 
${\rm CFT}(A_1)$ to a disorder correlator of ${\rm CFT}(A_2)$ and vice versa.

As an illustration, take $A_1 \,{=}\, {\bf 1}$ (so that ${\rm CFT}(A_1)$ 
is the Cardy case) and let $A_2$ be a Schellekens algebra. 
Set $B \,{=}\, A_2$, which becomes an $A_1$-$A_2$-bimodule
by taking the trivial action of $A_1$ for the left action
and the product of $A_2$ for the right action. Then
$B{\otimes_{\!A_2}}B^\vee \,{\cong}\, A_2$ is in
$\calp ic({}_{A_1\!}{\mathcal C}_{\!A_1}) \,{\equiv}\,
\calp ic(\calc)$ by definition, and one can also show that
$B^\vee{\otimes_{\!A_1}}B \,{=}\, B^\vee{\otimes}B$ is in 
$\calp ic({}_{A_2}\!{\,}{\mathcal C}_{\!A_2})$. Thus we can conclude that 
the correlators of any simple current CFT are related to the
correlators in the corresponding Cardy case by an order-disorder duality.


\section{Conclusions}

We have developed a rigorous algebraic approach to correlation functions in
rational conformal field theory. One aspect of this approach we have not
emphasized in this contribution is its computational power. Indeed, in the
construction of a rational conformal field theory for known chiral data,
there is only a single non-linear constraint to be solved: associativity
of the multiplication of the Frobenius algebra. The rest of the algorithm
is linear and allows, in the end, to find model-independent expressions for 
interesting CFT quantities like OPE coefficients \cite{fuRsn}
or coefficients of partition functions \cite{fuRs4,fuRs8}
in terms of invariants of knots and links in three-manifolds.

The algebraic approach we have presented in this contribution allows to
translate old physical problems to standard problems in algebra and 
representation theory. We conclude by summarizing this in the following
table:

\vskip .8em

\begin{center}
\begin{tabular}{lll}
Physical problem& $\longrightarrow$ &  Algebraic problem\\[.4em]
\hline \\[-.5em]
Classification of  CFTs & $\longrightarrow$ & Morita classes of algebras in $\calc$ 
\\
$\quad$~\,based on chiral data $\calc$ & & {\small
   (Many examples from Picard group of $\calc$)} 
\\[.4em]
Classification of &$\longrightarrow$ & Classification of
\\[.2em]
$\quad$\begin{tabular}{l} boundary conditions \\[.1em] defects
  \end{tabular} $\!\!\!\!\!\!$ \raisebox{-1.9pt}{\Huge \}}
  & & $\quad$\raisebox{-1.9pt}{\Huge \{}$\!\!\!\!$ \begin{tabular}{l}
  left modules \\[.1em] bimodules \end{tabular} 
\\[1.1em]
Internal symmetries & $\longrightarrow$ & Picard group of \calcaa\ 
\\[.4em]
Dualities & $\longrightarrow$ & Duality defects 
\\[.4em]
Deformation of CFTs & $\longrightarrow$ & Deformation of algebras 
{\small (and of categories)}
\\[.7em]
\end{tabular}
\end{center}

 \def\hy            {$\mbox{-\hspace{-.66 mm}-}$\linebreak[0]}
 \newcommand\wb{\,\linebreak[0]} \def\wB {$\,$\wb}
 \newcommand\Bi[2]    {\bibitem{#1}}
 \newcommand\JO[6]    {{\em #6}, {#1} {\bf #2} ({#3}), {#4--#5} }
 \newcommand\J[7]     {{\em #7}, {#1} {\bf #2} ({#3}), {#4--#5} {{\tt [#6]}}}
 \newcommand\JX[6]    {{\em #6}, {#1} {\bf #2} ({#3}), {#4} {{\tt [#5]}}}
 \newcommand\Prep[2]  {{\em #2}, pre\-print {\tt #1}}
 \newcommand\Pret[2]  {{\em #2}, pre\-print {#1}}
 \newcommand\BOOK[4]  {{\em #1\/} ({#2}, {#3} {#4})}
 \newcommand\inBO[7]  {{\em #7}, in:\ {\em #1}, {#2}\ ({#3}, {#4} {#5}), p.\ {#6}}
 \def\jf    {J.\ Fuchs}
 \def\adma  {Adv.\wb Math.}
 \def\anop  {Ann.\wb Phys.}
 \def\aspm  {Adv.\wb Stu\-dies\wB in\wB Pure\wB Math.}
 \def\atmp  {Adv.\wb Theor.\wb Math.\wb Phys.}
 \def\coma  {Con\-temp.\wb Math.}
 \def\comp  {Com\-mun.\wb Math.\wb Phys.}
 \def\cpma  {Com\-pos.\wb Math.}
 \def\duke  {Duke\wB Math.\wb J.}
 \def\fiic  {Fields\wB Institute\wB Commun.}
 \def\gafa  {Geom.\wB and\wB Funct.\wb Anal.}
 \def\ijmp  {Int.\wb J.\wb Mod.\wb Phys.\ A}
 \def\ijmb  {Int.\wb J.\wb Mod.\wb Phys.\ B}
 \def\jgap  {J.\wb Geom.\wB and\wB Phys.}
 \def\joag  {J.\wB Al\-ge\-bra\-ic\wB Geom.}
 \def\joal  {J.\wB Al\-ge\-bra}
 \def\jomp  {J.\wb Math.\wb Phys.}
 \def\jpaa  {J.\wB Pure\wB Appl.\wb Alg.}
 \def\maan  {Math.\wb Annal.}
 \def\nuci  {Nuovo\wB Cim.}
 \def\nupb  {Nucl.\wb Phys.\ B}
 \def\phlb  {Phys.\wb Lett.\ B}
 \def\phrl  {Phys.\wb Rev.\wb Lett.}
 \def\pspm  {Proc.\wb Symp.\wB Pure\wB Math.}
 \def\slnp  {Sprin\-ger\wB Lecture\wB Notes\wB in\wB Physics}
 \def\taia  {Topology \wB and\wB its\wB Appl.}
 \def\trgr  {Trans\-form.\wB Groups}
 \newcommand\mbop[2] {\inBO{The Mathematical Beauty of Physics}
            {J.M.\ Drouffe and J.-B.\ Zuber, eds.} \WS\Si{1997} {{#1}}{{#2}} }
 \newcommand\geap[2] {\inBO{Physics and Geometry} {J.E.\ Andersen, H.\
            Pedersen, and A.\ Swann, eds.} \MD\NY{1997} {{#1}}{{#2}} }
   \def\AMS    {{American Mathematical Society}}
   \def\AW     {{Addi\-son\hy Wes\-ley}}
   \def\Be     {{Berlin}}
   \def\Ca     {{Cambridge}}
   \def\CUP    {{Cambridge University Press}}
   \def\MD     {{Marcel Dekker}}
   \def\NY     {{New York}} 
   \def\OUP    {{Oxford University Press}}
   \def\PR     {{Providence}}
   \def\SV     {{Sprin\-ger Ver\-lag}}
   \def\WS     {{World Scientific}}
   \def\Si     {{Singapore}}

\bibliographystyle{amsalpha}

\end{document}